\documentclass[12pt,leqno]{article}
\usepackage{amsmath, amsthm, amsfonts, amssymb, color}
\usepackage{mathrsfs}
\theoremstyle{definition}

\newcommand{\scr}[1]{\mathscr #1}
\definecolor{wco}{rgb}{0.5,0.2,0.3}

\numberwithin{equation}{section} \theoremstyle{remark}

\newcommand{\ua}{\uparrow}

\title{{\bf Semigroup Properties for  the Second Fundamental Form}\footnote{Supported in
 part by NNSFC(10721091) and the 973-Project.}
}
\author{
{\bf Feng-Yu Wang}\\
\footnotesize{School of Mathematical Sci. and Lab. Math. Com. Sys.,
Beijing Normal
University, Beijing 100875, China}\\
\footnotesize{and}\\ \footnotesize{Department of Mathematics,
Swansea University, Singleton Park, SA2 8PP, UK}\\ \footnotesize{Email: wangfy@bnu.edu.cn;
F.Y.Wang@swansea.ac.uk}}
\begin{document}
\def\R{\mathbb R}  \def\ff{\frac} \def\ss{\sqrt} \def\BB{\mathbb
B}
\def\N{\mathbb N} \def\kk{\kappa} \def\m{{\bf m}}
\def\dd{\delta} \def\DD{\Delta} \def\vv{\varepsilon} \def\rr{\rho}
\def\<{\langle} \def\>{\rangle} \def\GG{\Gamma} \def\gg{\gamma}
  \def\nn{\nabla} \def\pp{\partial} \def\tt{\tilde}
\def\d{\text{\rm{d}}} \def\bb{\beta} \def\aa{\alpha} \def\D{\scr D}
\def\EE{\mathbb E} \def\si{\sigma} \def\ess{\text{\rm{ess}}}
\def\beg{\begin} \def\beq{\begin{equation}}  \def\F{\scr F}
\def\Ric{\text{\rm{Ric}}} \def\Hess{\text{\rm{Hess}}}\def\B{\scr B}
\def\e{\text{\rm{e}}} \def\ua{\underline a} \def\OO{\Omega}  \def\oo{\omega}     \def\tt{\tilde} \def\Ric{\text{\rm{Ric}}}
\def\cut{\text{\rm{cut}}} \def\P{\mathbb P} \def\ifn{I_n(f^{\bigotimes n})}
\def\C{\scr C}      \def\aaa{\mathbf{r}}     \def\r{r}
\def\gap{\text{\rm{gap}}} \def\prr{\pi_{{\bf m},\varrho}}  \def\r{\mathbf r}
\def\Z{\mathbb Z} \def\vrr{\varrho} \def\ll{\lambda}
\def\L{\scr L}\def\Tt{\tt} \def\TT{\tt}\def\II{\mathbb I}
\def\i{{\rm i}}\def\Sect{{\rm Sect}}\def\E{\mathbb E}

\maketitle
\begin{abstract}  Let $M$ be a compact Riemannian manifold with boundary
$\pp M$  and $L=
\DD+Z$ for a $C^1$-vector field $Z$ on $M$. Several equivalent statements, including the gradient and
Poincar\'e/log-Sobolev type inequalities of the Neumann semigroup generated by $L$, are presented for
lower bound conditions on the curvature of $L$ and
the second fundamental form of $\pp M$. The main result not only  generalizes the corresponding known
 ones on manifolds without boundary, but also clarifies the role of the second fundamental form in the
  analysis of the Neumann semigroup. Moreover, the L\'evy-Gromov isoperimetric inequality is also studied
  on manifolds with boundary. \end{abstract} \noindent
 \vskip 2cm

\section{Introduction}
The main purpose of this paper is to find out equivalent properties of the Neumann
semigroup on   manifolds with boundary  for  lower bounds  of  the second fundamental form of the boundary.
To explain the main idea of the study, let us  briefly  recall some equivalent semigroup properties for curvature lower bounds  on manifolds without boundary.

Let $M$ be a   connected complete Riemannian manifold
without boundary  and let $L=\DD+ Z$ for some $C^1$-vector field $Z$ on $M$.
Let $P_t$ be the diffusion semigroup generated by $L$, which is unique and
Markovian if the curvature of $L$ is bounded below, namely (see \cite{BE}),

\beq\label{C} \Ric-\nn Z  \ge -K\end{equation} holds on $M$ for some
constant $K\in \R$.
The following is a collection of known equivalent statements for
(\ref{C}), where the first two ones on gradient estimates are
classical in geometry (see e.g. \cite{B1, DL, DS, EL}), and the
remainder follows from Propositions 2.1 and 2.6 in \cite{B2}
(see also \cite{Le}):

\begin{enumerate}  \item[({\rm i})] $\ |\nn P_t f|^2\le \e^{2Kt}P_t|\nn
f|^2,\ \ t\ge 0,\ f\in C_b^1(M);$\item[({\rm ii})] $\ |\nn P_t f|\le
\e^{Kt}P_t|\nn f|,\ \ t\ge 0,\ f\in C_b^1(M);$\item[({\rm iii})] $\ P_t
f^2-(P_t f)^2\le\dfrac{\e^{2Kt}-1}{K}P_t|\nn f|^2,\ \ t\ge 0,\
f\in C_b^1(M);$ \item[({\rm iv})] $\ P_t f^2- (P_tf)^2 \ge
\dfrac{1-\e^{-2Kt}}{K} |\nn P_t f|^2,\ \ t\ge 0,\ f\in C_b^1(M);$
\item[({\rm v})] $\ P_t(f^2\log f^2)- (P_tf^2)\log (P_tf^2)\le
\dfrac{2(\e^{2Kt}-1)}K P_t|\nn f|^2,\ \ t\ge 0,\ f\in C_b^1(M);$
\item[({\rm vi})] $\ (P_t f)\{ P_t(f\log f)- (P_tf)\log
(P_tf)\}\ge \dfrac{1-\e^{-2Kt}}{2K}|\nn P_t f|^2,\ \ t\ge 0,\ f\in
C_b^1(M),\ f\ge 0.$\end{enumerate}

These equivalent statements for the curvature condition are crucial in the study
of heat semigroups and functional inequalities on manifolds. For the case that $M$
has a convex boundary, these equivalences are also true for $P_t$ the Neumann
semigroup (see \cite{W04} for  one more equivalent statement on Harnack
inequality). The question is now can we extend this result to manifolds with non-convex boundary,
and furthermore  describe  the  second fundamental
using semigroup properties?

So, from now on we assume that $M$ has a boundary
$\pp M$. Let $N$ be the inward unit normal vector field on $\pp M$. Then the second fundamental form
is a two-tensor on $T\pp M$, the tangent space of $\pp M$, defined by

$$\II(X,Y)= -\<\nn_X N, Y\>,\ \ \ X,Y\in T\pp M.$$
If $\II\ge 0$(i.e. $\II(X,X)\ge 0$ for $X\in T\pp M$), then $\pp M$ (or $M$) is called
convex. In general, we intend to study the lower bound condition of $\II$; namely,
$ \II\ge -\si$ on $\pp M$ for some $\si\in \R$.

For $x\in M$, let $\E^x$ be the expectation taken for the reflecting $L$-diffusion
process $X_t$ starting from $x$. So, for a bounded measurable  functional $\Phi$ of $X$,

$$\E\Phi:\ x\mapsto \E^x \Phi$$ is a function on $M$. Moreover, let $l_t$ be
the local time of $X_t$ on $\pp M$.
According to \cite[Theorem 5.1]{Hsu},
(\ref{C}) and $\II\ge -\si$  imply

\beq\label{G} |\nn P_t f|\le \e^{Kt} \E \big[|\nn f|(X_t)|\e^{\si l_t}\big],\ \ \ t>0, f\in C^1(M).
\end{equation}   To see that (\ref{G}) is indeed equivalent to  (\ref{C}) and $\II\ge -\si$,
we shall make use of the following formula for the
 second fundamental form established  recently by the author in  \cite{W09}:
 for any $f\in C^\infty(M)$ satisfying the Neumann condition $Nf|_{\pp M}=0$,

 \beq\label{W09} \II(\nn f,\nn f)= \ff{\ss \pi|\nn f|^2} 2 \lim_{t\to 0} \ff 1 {\ss t}
 \log \ff{(P_t|\nn f|^p)^{1/p}}{|\nn P_t f|}\end{equation} holds on $\pp M$ for
 any $p\in [1,\infty).$  With help of this result and stochastic analysis on the reflecting diffusion process,
 we are able to prove  the following   main result of the paper.

 \beg{thm}\label{T1.1} Let $M$ be a compact Riemannian manifold with boundary and let
 $P_t$ be the Neumann semigroup generated by $L=\DD+Z$. Then for any constants
 $K,\si\in \R$, the following
 statements are equivalent to each other: \begin{enumerate}
 \item[$(1)$]  $\Ric-\nn Z\ge -K$ on $M$ and $\II\ge -\si$ on $\pp M$;  \item[ $(2)$] $ (\ref{G})$ holds; \item[$(3)$] $|\nn P_t f|^2\le
\e^{2Kt}(P_t|\nn f|^2)\E\e^{2\si l_t},\ \ t\ge 0,\ f\in C^1(M);$ \item[ $(4)$] $ P_t
(f^2\log f^2)-(P_t f^2)\log P_tf^2\le 4\E\big[|\nn f|^2(X_t) \int_0^t \e^{2\si(l_t-l_{t-s})+2Ks}\d s\big],
 \newline\ \ \ \ \  \ t\ge 0,\
f\in C^1(M);$ \item[$(5)$] $P_t f^2- (P_tf)^2 \le 2\E\big[|\nn f|^2(X_t)
\int_0^t \e^{2\si(l_t-l_{t-s})+2Ks}\d s\big],  \ t\ge 0,\
f\in C^1(M);$\item[$(6)$] $\ |\nn P_t f|^2 \le \Big(\dfrac{2K}{1-\e^{-2Kt}}\Big)^2
\big(P_t(f\log f)- (P_tf)\log P_t f\big)\E \big[f(X_t)\int_0^t \e^{2\si l_s-2Ks}\d s\big],
 \newline\ \ \ \ \  \ t> 0,\  f\ge 0,\  f\in C^1(M);$ \item[$(7)$] $\ |\nn P_t f|^2
 \le \dfrac {2K^2}{(1-\e^{-2Kt})^2}
\big(P_t f^2-(P_tf)^2\big) \E\int_0^t \e^{2\si l_s-2Ks}\d s,\ \ t\ge 0,\ f\in C^1(M).$
\end{enumerate} \end{thm}

Theorem \ref{T1.1}  can be extended to a class of non-compact manifolds with boundary
such that the local times $l_t$ is exponentially integrable.  According to \cite{W09b}
the later is true provided $\II$ is bounded, the sectional curvature around $\pp M$ is bounded above,
the drift $Z$ is bounded around $\pp M$, and the injectivity radius
of the boundary is positive. To avoid technical complications, here we simply consider the
compact case.

In the next section, we shall provide a result on gradient estimate and non-constant lower bounds of
curvature and second fundamental form, which implies the equivalences
among (1), (2) and (3) as a special case. Then we present a complete proof for
the remainder of Theorem \ref{T1.1} in Section 3. As mentioned above, for manifolds without boundary
or with a convex boundary an equivalent Harnack inequality for
the curvature condition has been presented in \cite{W04}.  Due to unboundedness of the local time
 which causes an essential difficulty in the study of Harnack inequality,  the corresponding result
 for lower bound conditions
of the curvature and  the second fundamental form is  still open. Finally, as an extension to a result in
\cite{BL} where manifolds without boundary is considered, the L\'evy-Gromov isoperimetric inequality is derived in
Section 4 for  manifolds with boundary.

\section{Gradient estimate}

Let $K_1, K_2\in C(M)$ be such that

\beq\label{2.1} \Ric -\nn Z\ge -K_1\ \text{on}\ M,\ \ \II\ge -K_2\ \text{on}\ \pp M.
\end{equation} According to \cite[Theorem 5.1]{Hsu} this condition implies

\beq\label{G2} |\nn P_t f|\le \E\big[|\nn f|(X_t)\e^{\int_0^t K_1(X_s)\d s+\int_0^t K_2(X_s)
\d l_s}\big],\ \ \ t\ge 0, f \in C^1(M).\end{equation}
The main purpose of this section  to prove that these two statements are indeed
equivalent to each other.  To prove that (\ref{G2}) implies (\ref{2.1}), we need
the following results collected from \cite[Proof of Lemma 2.1]{W05} and
 \cite[Theorem 2.1, Lemma 2.2,  Proposition A.2]{W09} respectively:

 \begin{enumerate} \item[(I)] For any $\ll>0, \ \E\e^{\ll l_t}<\infty$.
 \item[(II)] For $X_0=x\in \pp M,\ \limsup_{t\to 0} \ff 1 t |\E l_t -2\ss{t/\pi}|<\infty.$
 \item[(III)] For $X_0=x\in \pp M$, there exists a constant $c>0$
 such that $\E l_t^2\le c t,\ \ t\in [0,1].$
 \item[(IV)] Let $\rr$ be the Riemannian distance. For $\dd>0$ and
 $X_0=x\in M\setminus \pp M$ such that $\rr  (x, \pp M)\ge \dd$, the stopping time
 $\tau_\dd:= \inf\{t>0: \rr(X_t,x)\ge \dd\}$ satisfies $\P(\tau_\dd\le t)\le c\exp[-\dd^2/(16 t)]$ for
 some constant $c>0$ and all $t>0.$ \end{enumerate}

\beg{thm}\label{T2.1} $(\ref{2.1})$, $(\ref{G2})$ and the following inequality are equivalent to each other:

\beq\label{G3}  |\nn P_t f|^2\le (P_t|\nn f|^2)\E\big[\e^{2\int_0^t K_1(X_s)\d s
+2\int_0^t K_2(X_s)
\d l_s}\big],\ \ \ t\ge 0, f \in C^1(M).\end{equation}\end{thm}

\beg{proof} Since by \cite{Hsu} (\ref{2.1}) implies (\ref{G2}) which is stronger than
(\ref{G3}) due to the Schwartz inequality, it remains to deduce (\ref{2.1}) from
(\ref{G3}).

(a) Proof of $\Ric-\nn Z\ge -K_1$. It suffices to prove at points in the interior. Let
$X_0=x\in M\setminus \pp M.$ For any $\vv>0$ there exists  $\dd>0$ such that

\beq\label{2.2} \bar B(x,\dd)\subset M\setminus \pp M,\ \ \ \sup_{y\in \bar B(x,\dd)}
|K_1(y)-K_1(x)|\le \vv,\end{equation} where $\bar B(x,\dd)$ is the closed
geodesic ball at $x$ with radius $\dd$. Since $l_t=0$ for $t\le \tau_\dd$, by (\ref{G3}), (I) and (IV) we have

\beg{equation*}\beg{split} |\nn P_t f|^2(x) &\le (P_t |\nn f|^2(x)) \E \e^{2
\int_0^t K_1(X_s)\d s + 2 \int_0^t K_2(X_s)\d l_s}\\
&\le (P_t |\nn f|^2(x)) \Big\{ \e^{2t(K_1(x)+\vv)}\P(\tau_\dd\ge t) +
\ss{\P(\tau_\dd<t) \E\e^{4t \|K_1\|_\infty + 4 \|K_2\|_\infty l_t}}\Big\}\\
&\le (P_t |\nn f|^2(x)) \e^{2t(K_1(x)+\vv)} + C \e^{-\ll/t},\ \ \ t\in (0,1]
\end{split}\end{equation*} for some constants $C,\ll>0.$ This implies

\beq\label{2.3} \limsup_{t\to 0} \ff{|\nn P_t f|^2(x) -|\nn f|^2(x)}t
\le \limsup_{t\to 0} \ff {\e^{2t(K_1(x)+\vv)}P_t |\nn f|^2(x) -|\nn f|^2(x)} t.\end{equation}
 Now, let $f\in C^\infty(M)$ with $Nf|_{\pp M}=0$, we have

 $$P_tf= f+ \int_0^t P_s Lf \d s,\ \ \ t\ge 0.$$ Then

 \beg{equation}\label{2.0}\beg{split} &\limsup_{t\to 0} \ff{|\nn P_t f|^2(x) -|\nn f|^2(x)}t \\
 = & \lim_{t\to 0} \ff 1 t \bigg\{\bigg|\int_0^t \nn P_s Lf\d s\bigg|^2 + 2 \int_0^t \<\nn f, \nn
 P_s Lf\>\d s\bigg\}(x).\end{split}\end{equation} Moreover, according to the last display
 in the proof of \cite[Theorem 5.1]{Hsu} (the initial data $u_0\in O_x(M)$ was missed in
 the right hand side therein),

 $$\nn P_t Lf= u_0 \E\big[M_tu_t^{-1}\nn L f(X_t)\big],$$ where $u_t$ is the horizontal
 lift of $X_t$ on the frame bundle $O(M)$, and $M_t$ is a $d\times d$-matrices valued
 right continuous process satisfying $M_0=I$ and  (see \cite[Corollary 3.6]{Hsu})

 $$\|M_t\|\le \exp\big[\|K_1\|_\infty t+\|K_2\|_\infty l_t\big].$$ So,
 due to (I), $|\nn P_\cdot Lf|$
 is bounded on $[0,1]\times M$ and $\nn P_s Lf\to \nn Lf$ as $s\to 0$. Combining
 this with (\ref{2.0}) we obtain

\beq\label{2.4}  \limsup_{t\to 0} \ff{|\nn P_t f|^2(x) -|\nn f|^2(x)}t= 2\<\nn f, \nn Lf\>(x).
\end{equation}
On the other hand, applying the It\^o formula to $|\nn f|^2(X_t)$ we have

\beq\label{2.5}\beg{split} P_t |\nn f|^2(x)& =|\nn f|^2(x) +\int_0^t P_s L|\nn f|^2(x)\d s +\E \int_0^t
N|\nn f|^2(X_s)\d l_s\\
&\le |\nn f|^2(x) +\int_0^t P_s L|\nn f|^2(x)\d s +\|\nn |\nn f|^2\|_\infty
 \E l_t.\end{split}\end{equation}
 Since $l_t=0$ for $t\le \tau_\dd$, by (III) and (IV) we have

 $$\E l_t \le \ss{(\E l_t^2)\P(\tau_\dd\le t)}\le c_1\e^{-\ll/t},\ \ \ t\in (0,1]$$ for some constants $c_1,\ll>0.$
  So, it follows from (\ref{2.5}) that

 $$\limsup_{t\to 0} \ff{P_t|\nn f|^2(x)-|\nn f|^2(x)} t \le L|\nn f|^2(x).$$ Combining this
 with (\ref{2.3}) and (\ref{2.4}), we arrive at

 $$\ff 1 2 L|\nn f|^2(x)- \<\nn f, \nn Lf\>(x)\ge -(K_1(x)+\vv),\ \ \ f\in C^\infty(M), Nf|_{\pp M}=0.$$
 According to the Bochner-Weitzenb\"ock formula, this is equivalent to
 $(\Ric-\nn Z)(x)\ge -(K_1(x)+\vv).$ Therefore, $\Ric-\nn Z\ge -K_1$ holds on $M$ by the
 arbitrariness of $x\in M\setminus \pp M$ and $\vv>0.$

 (b) Proof of $\II\ge -K_2$. Let $X_0=x\in \pp M.$ For any $f\in C^\infty(M)$ with
 $Nf|_{\pp M}=0$, (\ref{G3}) implies that

 \beq\label{2.6} |\nn P_t f|^2 (x)\le \e^{C_1t} (P_t |\nn f|^2(x))
 \E\e^{2\int_0^t K_2(X_s)\d l_s},
 \end{equation} where $C_1= 2\|K_1\|_\infty.$  Let

 $$\vv_t= 2\sup_{s\in [0,t]} |K_2(X_s)-K_2(x)|.$$ By the continuity of the reflecting diffusion
 process we have $\vv_t\downarrow 0$ as $t\downarrow 0.$ Since there exists $c_0>0$
 such  that for any $r\ge 0$ one has
 $\e^r\le 1 + r +c_0 r^{3/2} \e^r,$ we obtain

\beq\label{2.7}\log \E \e^{2\int_0^t K_2(X_s)\d l_s} \le \log\big\{1 + 2 K_2(x) \E l_t + \E(\vv_t l_t)
 + C_2 \E(l_t^{3/2}\e^{C_2 l_t})\big\}\end{equation} for some constant $C_2>0.$  Moreover, by (I) and (III)
 we have

 $$ \E(l_t^{3/2} \e^{C_2 l_t})\le (\E l_t^2)^{3/4} (\E \e^{4C_2 l_t})^{1/4}\le C_3 t^{3/4},\ \ \
 t\in (0,1]$$ for some constant $C_3>0.$  Substituting this and (\ref{2.7}) into (\ref{2.6}),
 we arrive at

 $$\limsup_{t\to 0}\ff 1 {\ss t} \log \ff{|\nn P_t f|^2(x)}{P_t|\nn f|^2(x)}
 \le \limsup_{t\to 0} \ff{2K_2(x) \E l_t +\E(\vv_t l_t)}{\ss t}.$$ Since $\E \vv_t^2\to 0$ as $t\to 0$ and
 $\E l_t^2\le c t$ due to (III), this and (II) imply

  $$\limsup_{t\to 0}\ff 1 {\ss t} \log \ff{|\nn P_t f|^2(x)}{P_t|\nn f|^2(x)}\le \ff{4K_2(x)}{\ss\pi}.$$
  Combining this with (\ref{W09}) for $p=2$ we complete the proof.
 \end{proof}

 \section{Proof of Theorem \ref{T1.1}}

 Applying   Theorem \ref{T2.1} to $K_1=K$ and $K_2 =\si$ we conclude that
 (1), (2) and (3) are equivalent to each other. Noting that the log-Sobolev inequality (4) implies
 the Poincar\'e inequality (5) (see e.g.  \cite{DS}),
  it suffices to prove that $(2)\Rightarrow (4)$, $(5)\Rightarrow (1)$, and $(2) \Rightarrow (6)\Rightarrow (7)\Rightarrow (1),$
  where $``\Rightarrow$'' stands for $``$implies''.
  We shall complete the proof step by step.

  (a) (2) $\Rightarrow$ (4). By approximations we may assume that $f\in C^\infty(M)$ with $Nf|_{\pp M}=0.$
  In this case

  $$\ff{\d}{\d t} P_t f= L P_t f= P_t Lf .$$ So, for fixed $t>0$ it follows from (2) that

  \beq\label{3.1}\beg{split}  \ff{\d}{\d s} P_{t-s}\{(P_s f^2)\log P_s f^2\} &=
  -P_{t-s} \ff{|\nn P_s f^2|^2}{P_sf^2}
   \ge -4 \e^{2Ks}P_{t-s}\ff{(\E [f|\nn f|(X_s)\e^{\si l_s}])^2}{P_s f^2}\\
   &\ge -4 \e^{2Ks} P_{t-s} \E [ |\nn f|^2(X_s) \e^{2\si l_s}].\end{split}\end{equation}
   Next, by the Markov property, for $\F_s= \si(X_r: r\le s), s\ge 0, $ we have

   \beg{equation*}\beg{split} P_{t-s} (\E[|\nn f|^2(X_s)\e^{2\si l_s}])(x)&= \E^x \E^{X_{t-s}} [|\nn f|^2(X_s)
   \e^{2\si l_s}]\\
   &= \E^x [\E^x(\e^{2\si (l_t-l_{t-s})}|\nn f|^2(X_t)|\F_{t-s})]= \E^x [|\nn f|^2(X_t)\e^{2\si (l_t-l_{t-s})}].
   \end{split}\end{equation*}
   Combining this with (\ref{3.1}) we obtain

   $$\ff{\d}{\d s} P_{t-s} \{(P_sf^2)\log P_s f^2\} \ge - 4 \E \big[|\nn f|^2(X_t) \e^{
   2Ks+ 2\si (l_t-l_{t-s})}\big],\ \ \ s\in (0,t).$$ This implies (4) by integrating both sides
   with respect to $\d s$ from $0$ to $t$.

   (b1) (5) $\Rightarrow\ \Ric-\nn Z\ge -K.$ Let $X_0=x\in M\setminus \pp
   M$ and $f\in C^\infty(M)$ with $Nf|_{\pp M}=0$. By (5) we have

   \beq\label{3.2} P_t f^2-(P_tf)^2 \le 2 \E\bigg[|\nn f|^2(X_t)
   \int_0^t \e^{2Ks+ 2\si(l_t-l_{t-s})}\d s\bigg].\end{equation} Let
   $\dd>0$ and $\tau_\dd$ be as in the proof of Theorem
   \ref{T2.1}(a). Then

   \beg{equation*}\beg{split} &\E\bigg[|\nn f|^2(X_t)
   \int_0^t \e^{2Ks+ 2\si(l_t-l_{t-s})}\d s\bigg]\\
   &\le (P_t|\nn f|^2)\int_0^t \e^{2Ks}\d s + t \|\nn f\|_\infty
   \e^{2Kt} \E[\e^{2\si l_t}1_{\{\tau_\dd<t\}}]\\
   &\le \ff{\e^{2Kt}-1}{2K} P_t|\nn f|^2(x) +
   c\e^{-\ll/t},\ \ \ t\in (0,1]\end{split}\end{equation*} holds for some constants
   $c,\ll >0$ according to (IV). Combining this with (\ref{3.2}) we
   conclude that

   \beq\label{3.3} P_tf^2(x)-(P_tf)^2(x)\le \ff{\e^{2Kt}-1}K
   P_t|\nn f|^2(x) + 2c\e^{-\ll/t},\ \ \ t\in (0,1].\end{equation}
   Since $f\in C^\infty(M)$ with $Nf|_{\pp M=0}$, we have

 \beg{equation}\label{A1}\beg{split} P_t f^2-(P_tf)^2 &= f^2 +\int_0^tP_s
 Lf^2\d s -\bigg(f+\int_0^t P_s Lf\d s\bigg)^2\\
 &= \int_0^t (P_s Lf^2 -2fP_s Lf)\d s -\bigg(\int_0^tP_s Lf\d
 s\bigg)^2.\end{split}\end{equation}
 Moreover, by the continuity of $s\mapsto P_s Lf$, we have

\beq\label{A2} \bigg(\int_0^tP_s Lf\d s\bigg)^2= (Lf)^2 t^2
+\circ(t^2),\end{equation}
 where and in what follows, for a positive function $(0,1]\ni t\mapsto \xi_t$ the notion
 $\circ(\xi_t)$ stands for a variable such that
 $\circ(\xi_t)/\xi_t\to 0$ as $t\to 0$; while $\bigcirc(\xi_t)$
 satisfies that $\bigcirc(\xi_t)/\xi_t$ is bounded for $t\in (0,1].$
Moreover, since

\beg{equation*}\beg{split}  P_s Lf^2- 2 f P_s Lf= &L f^2- 2fLf +\int_0^s
 (P_r L^2 f^2-2f P_r
 L ^2 f)\d r\\
 & + \E\int_0^s (N Lf^2- 2f(x) N Lf)(X_r)\d l_r, \end{split}\end{equation*}  and due
 to (IV)

$$\bigg|\E\int_0^t \big\{N Lf^2- 2f(x) N Lf\big\}(X_r)\d l_r\bigg|\le c_1 \E
l_s \le c_2 \e^{-\ll/s},\ \ \ s\in (0,1]$$ holds for some constants $c_1,
c_2, \ll>0$, it follows from the continuity of $P_s $ in $s$ that

$$ \int_0^t (P_s Lf^2 -2fP_s Lf)\d s  = 2t |\nn f|^2 +\ff {t^2} 2 (L^2f^2-2fL^2f)+
\circ(t^2).$$ Combining this with (\ref{A1}) and (\ref{A2}) we
obtain

\beq\label{3.4} \beg{split} P_tf^2(x)-(P_tf)^2(x)&= 2t|\nn f|^2(x)
+\ff{t^2} 2
(L^2f^2- 2fL^2f)(x) -t^2(Lf)^2(x) +\circ(t^2)\\
&= 2t |\nn f|^2(x) +t^2(2\<\nn f, \nn Lf\>+L|\nn f|^2)(x) +\circ
(t^2).\end{split}\end{equation} Similarly,

\beg{equation*}\beg{split} P_t|\nn f|^2(x)&= |\nn f|^2(x)
+\int_0^tP_s L|\nn f|^2(x)\d s +\E
\int_0^t N|\nn f|^2(X_s)\d l_s\\
&=|\nn f|^2(x) +t L|\nn f|^2(x) +\circ(t).\end{split}\end{equation*}
Combining this with (\ref{3.3}) and (\ref{3.4}) we arrive at

\beg{equation*}\beg{split} &\ff 1 {t^2} \big\{t^2 (2\<\nn f, \nn
Lf\>+L|\nn f|^2)(x)+\circ(t^2)\big\}\\
& \le \ff{\e^{2Kt}-1}{Kt}L|\nn f|^2(x) +\circ(1) +\ff 1 t
\Big(\ff{\e^{2Kt}-1}{Kt}-2\Big)|\nn
f|^2(x).\end{split}\end{equation*} Letting $t\to 0$ we obtain

$$L|\nn f|^2(x)-2\<\nn f,\nn Lf\>(x) \ge -2K|\nn f|^2(x),$$ which
implies $(\Ric-\nn Z)(x)\ge -K$ by the Bochner-Weitzenb\"ock
formula.

(b2) (5) $\Rightarrow\ \II\ge -\si.$ Let $X_0=x\in \pp M$ and $f\in C^\infty(M)$ with $Nf|_{\pp M}=0$.
Noting that $L f^2-2fLf=2|\nn f|^2$, by the It\^o formula  we have

\beg{equation}\label{3.5}\beg{split}& P_t f^2(x)-(P_t f)^2(x)= f^2
 +\int_0^t P_s Lf^2\d s-\bigg( f+ \int_0^t P_s Lf\d s\bigg)^2\\
 &= 2\int_0^t P_s |\nn f|^2(x)\d s + 2 \int_0^t [P_s(fLf)(x)- f(x) P_s Lf(x)]\d s+\bigcirc(t^2).\end{split}\end{equation}
 Since $Nf|_{\pp M}=0$ implies

 $$0=\<\nn f,\nn \<N,\nn f\>\>= \Hess_f(N,\nn f)-\II(\nn f,\nn f),$$
it follows that

 \beq\label{**}\II(\nn f,\nn f)= \Hess_f(N,\nn f)= \ff 1 2 N|\nn f|^2.\end{equation} So, by the It\^o formula, (II) and (III) yield

 \beg{equation}\label{3.6}\beg{split} P_s |\nn f|^2 (x)
&= |\nn f|^2(x) +\int_0^s P_r L|\nn f|^2(x)\d r
  +\E \int_0^s N|\nn f|^2(X_r)\d l_r\\
  &= |\nn f|^2(x) +\bigcirc(s) + 2 \E\int_0^s \II(\nn f,\nn f)(X_r)\d l_r\\
 &= |\nn f|^2(x) + \ff {4\ss s}{\ss\pi} \II(\nn f,\nn f)(x)  + \circ(s^{1/2}).\end{split}\end{equation}
Moreover,  since $(fNLf)(X_r)-f(x)(NLf)(X_r)$ is bounded and goes to zero as $r\to 0$, it follows from (III) that

 $$ 2 \E \int_0^t \d s \int_0^s [(fNf)(X_r) -f(x)(NLf)(X_r)]\d l_r= \circ(t^{3/2}).$$ So, by the I\^o formula

 \beg{equation*}\beg{split} & 2 \int_0^t [P_s(fLf)(x)- f(x) P_s Lf(x)]\d s\\
 &=  2\int_0^t \d s \int_0^s [P_r L(fLf)(x)- f(x) P_r L^2 f(x)]
 \d r\\
 &\qquad + 2 \E \int_0^t \d s \int_0^s [(fNLf)(X_r) -f(x)(NLf)(X_r)]\d l_r = \circ(t^{3/2}).\end{split}\end{equation*}
   Combining this with (\ref{3.5}) and  (\ref{3.6}) we arrive at

 \beq\label{3.8}\beg{split}
 &\lim_{t\to 0} \ff 1 {t\ss t} \big(P_tf^2(x)-(P_tf)^2(x) - 2 t |\nn f|^2(x)\big)\\
  &= \ff 8{\ss\pi} \II(\nn f,\nn f)(x)\lim_{t\to 0} \ff 1 {t\ss t} \int_0^t \ss s\, \d s
 =\ff{16}{3\ss\pi} \II(\nn f,\nn f)(x).\end{split}\end{equation} On the other hand,
 by the It\^o formula for $|\nn f|^2(X_t)$, it follows from (\ref{**}) and (II) that

\beg{equation}\label{3.*}\beg{split} A_t& := \ff 1 {t\ss t }\E \bigg\{|\nn f|^2(X_t) \int_0^t
 \e^{2Ks+ 2 \si(l_t-l_{t-s})}\d s- t |\nn f|^2(x)\bigg\}\\
&=\ff 1 {\ss t} \big(\E |\nn f|^2(X_t) -|\nn f|^2(x)\big) +  \E\bigg\{ \ff{|\nn f|^2(X_t)}
{t\ss t}\int_0^t \big(\e^{2Ks+2\si (l_t-l_{t-s})}-1\big)\d s\bigg\}\\
&= \ff 1 {\ss t}   \bigg\{\int_0^t P_sL|\nn f|^2(x)\d
s+ \E \int_0^t N|\nn f|^2(X_s)\d l_s\bigg\}\\
&\qquad  +\E\bigg\{ \ff{|\nn f|^2(X_t)}
{t\ss t}\int_0^t \big(\e^{2Ks+2\si (l_t-l_{t-s})}-1\big)\d s\bigg\}\\
&= \ff 4 {\ss\pi} \II(\nn f,\nn f)(x)+\circ(1) +\E\bigg\{ \ff{|\nn f|^2(X_t)}
{t\ss t}\int_0^t \big(\e^{2Ks+2\si (l_t-l_{t-s})}-1\big)\d s\bigg\}.\end{split}
\end{equation} Since by (I) and (III)

\beg{equation*}\beg{split} & \bigg|\E \Big[\big(|\nn f|^2(X_t)-|\nn f|^2(x)\big)\int_0^t \big(\e^{2Ks+2\si (l_t-l_{t-s})}-1\big)\d s\Big]\bigg|\\
&\le t \Big\{\E\big(|\nn f|^2(X_t)-|\nn f|^2(x)\big)^2\Big\}^{1/2} \Big\{\E \big(\e^{2Kt+2\si l_t}-1\big)^2\Big\}^{1/2}\\
&= \circ(t) \cdot \big(\E[4\si^2 l_t^2] +\circ(t)\big)= \circ (t^2),\end{split}\end{equation*} it follows from (I) and (II) that

\beg{equation*}\beg{split}&\E \bigg[|\nn f|^2(X_t)\int_0^t  \big(\e^{2Ks+2\si (l_t-l_{t-s})}-1\big)\d s\bigg]\\
&= \circ(t^2) + |\nn f|^2(x) \E \int_0^t  \big(\e^{2Ks+2\si (l_t-l_{t-s})}-1\big)\d s\\
&=\circ(t^{3/2}) +\ff{4\si |\nn f|^2(x) }{\ss\pi} \int_0^t \big(\ss t-\ss{t-s}\big)\d s\\
&= \ff{4\si t\ss t}{3\ss\pi} |\nn f|^2(x)+ \circ (t^{3/2}).\end{split}\end{equation*}
Combining this with (\ref{3.*})  we arrive at

$$ A_t \le \circ(1) + \ff 4{\ss\pi} \II(\nn f,\nn f)(x) + \ff{ 4\si}{3\ss\pi}|\nn f|^2(x).$$
So,  (\ref{3.8})
and   (5) imply that

$$\ff{16}{3\ss\pi} \II(\nn f,\nn f)(x)\le \limsup_{t\to 0}2A_t \le \ff{8}{\ss\pi} \II(\nn f,\nn f)(x)
+\ff{8\si}{3\ss\pi} |\nn f|^2(x).$$ Therefore,
$\II(\nn f,\nn f)(x)\ge -\si |\nn f|^2(x).$

(c) (2) $\Rightarrow$ (6).  Let $f\ge 0$ be smooth satisfying the Neumann boundary
condition. We have

$$\ff{\d}{\d s} P_s \big\{(P_{t-s}f)\log P_{t-s} f\big\} = P_s \ff{|\nn P_{t-s}f|^2}{P_{t-s}f}.$$
This implies

\beq\label{B1} P_t(f\log f) -(P_t f)\log P_t f=\int_0^t P_s \ff{|\nn P_{t-s}f|^2}{P_{t-s} f}\d s.
\end{equation} On the other hand, by (2) and applying the Schwartz inequality
to the probability measure $\ff{2K}{1-\exp[-2Kt]}\e^{-2Ks}\d s$ on $[0,t],$
we obtain

\beg{equation*}\beg{split} |\nn P_t f|^2 &= \bigg\{\ff{2K}{1-\e^{-2Kt}}
\int_0^t |\nn P_s(P_{t-s}f)|\e^{-2Ks}\d s\bigg\}^2\\
&\le  \bigg\{\ff{2K}{1-\e^{-2Kt}}
\int_0^t E\big[|\nn P_{t-s}f|(X_s)\e^{\si l_s- Ks}\big]\d s\bigg\}^2\\
&\le \Big(\ff{2K}{1-\e^{-2Kt}}\Big)^2 \bigg(\E\int_0^t \ff{|\nn P_{t-s}f|^2}{P_{t-s}f}(X_s)\d s\bigg)
\int_0^t \E \big[P_{t-s}f(X_s) \e^{2\si l_s-2Ks}\big]\d s\\
&=  \Big(\ff{2K}{1-\e^{-2Kt}}\Big)^2 \bigg(\int_0^tP_s \ff{|\nn P_{t-s}f|^2}{P_{t-s}f}\d s\bigg)\int_0^t
 \E \big[P_{t-s}f(X_s) \e^{2\si l_s-2Ks}\big]\d s.\end{split}\end{equation*} Combining this with (\ref{B1}) and noting that
 the Markov property implies

\beg{equation*}\beg{split} \E [P_{t-s}f(X_s)\e^{2\si l_s}] &= \E [(\E^{X_s}f(X_{t-s}))
\e^{2\si l_s}]= \E[\e^{2\si l_s} \E(f(X_t)|\F_s)]\\
&= \E[\E(f(X_t)\e^{2\si l_s}|\F_s)] = \E[f(X_t)\e^{2\si l_s}],\end{split}\end{equation*} we obtain (6).

(d) (6) $\Rightarrow$ (7). The proof is similar to the classical one for the log-Sobolev inequality to imply
the Poincar\'e inequality. Let $f\in C^\infty(M)$. SInce $M$ is compact,  $1+\vv f> 0$ for small $\vv>0.$ Applying (6)
 to $1+\vv f$ in place of $f$, we obtain

\beq\label{B2} \beg{split} |\nn P_t f|^2 \le &  \ff{2K}{\vv^2(1-\e^{-2Kt})}
\big\{P_t (1+\vv f)\log(1+\vv f)- (1+\vv P_t f)\log (1+\vv P_t f)\big\} \\
&\cdot \E \bigg\{(1+\vv f(X_t)) \int_0^t\e^{2\si  l_s-2Ks}\d s\bigg\}.\end{split}\end{equation}
Since by Taylor's expansion

$$P_t (1+\vv f)\log (1+\vv f) - (1+\vv P_t f)\log (1+\vv P_t f)=\ff {\vv^2} 2 \big(P_tf^2-
(P_t f)^2\big) + \circ(\vv^2),$$ letting $\vv\to 0$ in (\ref{B2}) we obtain (7).

(e1) (7) $\Rightarrow\ \Ric-\nn Z\ge -K.$ Let $X_0=x\in M\setminus\pp M$
and $f\in C^\infty(M)$ with $Nf|_{\pp M}=0$.  by (I) and (IV) we have

$$ \E \e^{2\si l_s} = 1 + \E[ \e^{2\si l_s}1_{\{\tau_\dd \le
s\}}]=1+\circ(s).$$ So,

$$\E \int_0^t \e^{2\si l_s-2Ks}\d s= \ff{1-\exp[-2Kt]}{2K} +\circ(t).$$ Combining this with
 (\ref{3.4}) and (7),
we conclude that,  at point $x$,

\beg{equation*} \beg{split} \ff{|\nn P_t f|^2-|\nn f|^2} t & \le \ff{K}{1-\e^{-2Kt}}
\big\{2|\nn f|^2+ t \big(2\<\nn f,\nn Lf\> +L|\nn f|^2\big)\big\}- \ff{|\nn f|^2} t +\circ(1)\\
&=\ff 1
t \Big(\ff{2Kt}{1-\e^{-2Kt}}-1\Big)|\nn f|^2 +    \ff{K t}{1-\e^{-2Kt}} \big(2\<\nn f, \nn Lf\> +L|\nn
f|^2\big) +\circ (1).\end{split}\end{equation*} Letting $t\to 0$   and using
(\ref{2.4}), we obtain

$$ 2\<\nn f,\nn Lf\>\le K|\nn f|^2+\<\nn f, \nn Lf\> +\ff 1 2 L|\nn
f|^2$$ at point $x$. This implies $\Ric-\nn Z\ge -K$ at this point
according to the Bochner-Weitzenb\"ock formula.

(e2) (7) $\Rightarrow\ \II\ge -\si.$ Let $X_0=x\in \pp M$ and $f\in
C^\infty(M)$ with $Nf|_{\pp M}=0.$  It follows from (\ref{3.8}), (7)
and (II) that at point $x$,

\beg{equation*}\beg{split} |\nn P_t f|^2 &\le \ff{2K^2}{(1-\e^{-2Kt})^2}
\Big(2t|\nn f|^2 + \ff{16t^{3/2}}{3\ss\pi} \II(\nn f,\nn
f)+\circ(t^{3/2})\Big)\Big(t+ \ff{8\si t^{3/2}}{3\ss\pi} +
\circ(t^{3/2})\Big)\\
&= \ff{4K^2t^2}{(1-\e^{-2Kt})^2} |\nn f|^2 +
\ff{4K^2t^{5/2}}{(1-\e^{-2Kt})^2}\Big(\ff 8 {3\ss\pi} \II(\nn
f,\nn f)+ \ff{8\si}{3\ss\pi} |\nn f|^2\Big) +
\circ(t^{1/2}).\end{split}\end{equation*}
Combining this with (\ref{2.4}) we deduce at point $x$ that

\beg{equation*}\beg{split} 0&= \lim_{t\to 0} \ff 1 {\ss t} \Big(
|\nn P_t f|^2-\ff{4K^2t^2}{(1-\e^{-2Kt})^2}|\nn f|^2\Big)\\
&\le \lim_{t\to 0} \ff{4K^2t^2}{(1-\e^{-2Kt})^2} \Big(\ff 8{3\ss\pi}
\II(\nn f,\nn f) +\ff {8\si}{3\ss\pi} |\nn f|^2\Big)\\
&= \ff 8 {3\ss\pi} \II(\nn f,\nn f) +\ff {8\si}{3\ss\pi} |\nn
f|^2.\end{split}\end{equation*} Therefore, $\II(\nn f,\nn f)(x)\ge
-\si |\nn f|^2(x).$

\section{L\'evy-Gromov isoperimetric inequality}
\def\U{\scr U}

As a dimension-free version of the classical L\'evy-Gromov isoperimetric
inequality, it is proved in \cite{BL} that if $M$ does not have boundary then for $V\in C^2(M)$ such that
$\Ric-\Hess_V\ge R>0$ the following inequality
\beq\label{4.1}\U(\mu(f))\le \int_M\ss{\U^2(f)+R^{-1} |\nn f|^2}\,\d\mu,\end{equation}
holds for any smooth function $f$ with values in $[0,1]$,
where $\mu(\d x):= C(V)^{-1}\e^{V(x)}\d x$ for $C(V)= \int_M
\e^{V(x)}\d x$ is a probability measure on $M$, and $\U= \varphi\circ\Phi^{-1}$ for $\Phi(r)= (2\pi)^{-1}\int_{-\infty}^r \e^{-s^2/2}\d s$
and $\varphi= \Phi'.$ Since $\U(0)=\U(1)=0,$ taking $f=1_A$ (by approximations) in (\ref{4.1}) for a smooth domain $A\subset M$, we obtain the isoperimetric inequality

\beq\label{4.1'}R\U(A)\le \mu_\pp (\pp A),\end{equation} where $\mu_\pp(\pp A)$ is the area of $\pp A$ induced by $\mu.$ This inequality is crucial
in the study of Gaussian type concentration of $\mu$ (see \cite{BL, Le}).  Obviously, (\ref{4.1}) follows from the following
semigroup inequality by letting $t\to\infty$:

\beq\label{4.2} \U(P_tf)\le P_t \ss{\U^2(f)+R^{-1}(1-\e^{-2Rt})|\nn f|^2}.\end{equation}
In this section we aim to extend (\ref{4.2}) to manifolds with boundary.

Now, let again $M$ be compact with boundary $\pp M$, and let $P_t$ be the Neumann semigroup generated by $L=\DD+Z$.
We shall prove an analogue of (\ref{4.2}) for the curvature and second fundamental condition in Theorem \ref{T1.1}(1).

\beg{thm}\label{T4.1} Let $\Ric-\nn Z\ge -K$ and $\II\ge -\si$ for some constants $K\in \R$ and $\si\ge 0.$ Then for any smooth function
$f$ with values in $[0,1]$,

\beq\label{4.3}  \U(P_tf)\le \E\ss{\U^2(f)(X_t) + |\nn f|^2(X_t)\ff{(\e^{2Kt}-1)\e^{2\si l_t}} K},\ \ \ t\ge 0.\end{equation}
If in particular $\pp M$ is convex $($i.e.  $\si=0)$, then

$$ \U(P_tf)\le P_t \ss{\U^2(f) + |\nn f|^2(X_t)\ff{\e^{2Kt}-1} K},\ \ \ t\ge 0.$$ If moreover $K<0$,
then $(\ref{4.1})$ and $(\ref{4.1'})$ hold for $R=-K>0.$\end{thm}

\beg{proof} It suffices to prove the first assertion. To this end, we shall use the following equivalent condition for $\Ric-\nn Z\ge -K$
(see e.g. the proof of \cite[(1.14)]{Le}):

\beq\label{C1} \GG_2( f,  f):= \ff 1 2 L|\nn f|^2 -\<\nn f, \nn Lf\>\ge -K|\nn f|^2 +\ff{|\nn|\nn f|^2|^2}{4|\nn f|^2}.\end{equation}
To prove (\ref{4.3}), we consider the process

$$\eta_s= \U^2(P_{t-s}f)(X_s) +|\nn P_{t-s}f|^2(X_s)\ff{(\e^{2Ks}-1)\e^{2\si l_s}}K,\ \ \ s\in [0,t].$$
To apply  the It\^o formula for $\eta_s$,  recall that $X_s$ solves the equation

$$\d X_s= \ss 2 \, u_s\circ \d B_s +N(X_s)\d l_s,$$ where $u_s$ is the horizontal lift of $X_s$ and $B_s$ is the Brownian motion
on $\R^d$ provided $M$ is $d$-dimensional. So,

\beg{equation*}\beg{split} &\d \eta_s = \ss 2 \Big\<2(\U \U' )(P_{t-s}f)(X_s)
+ \ff{(\e^{2Ks}-1)\e^{2\si l_s}}K\nn |\nn P_{t-s}f|^2(X_s), u_s\d B_s\Big\>\\
&\quad + \Big\{2({\U'}^2+\U\U'')(P_{t-s}f)|\nn P_{t-s}f|^2 + 2 \GG_2(P_{t-s}f, P_{t-s}f)\ff{(\e^{2Ks}-1)\e^{2\si l_s}}K\\
&\quad + 2 |\nn P_{t-s}f|^2
\e^{2Ks+2\si l_s}\Big\}(X_s)\d s + \ff{(\e^{2Ks}-1)\e^{2\si l_s}}K \big(N|\nn P_{t-s}f|^2 + 2\si|\nn P_{t-s}f|^2\big)(X_s)\d l_s.\end{split}\end{equation*}
Noting that $\U\U''=-1$ and $\si\ge 0$ so that $\e^{2\si l_s}\ge 1$, combining this with (\ref{**}), $\II\ge -\si$  and (\ref{C1}), we obtain

\beg{equation*}\beg{split} \d \eta_s &\ge \ss 2 \Big\<2(\U \U' )(P_{t-s}f)(X_s)
+ \ff{(\e^{2Ks}-1)\e^{2\si l_s}}K\nn |\nn P_{t-s}f|^2(X_s), u_s\d B_s\Big\>\\
&\quad +\Big\{2{ \U'}^2(P_{t-s}f) |\nn P_{t-s}f|^2 +\ff{(\e^{2Ks}-1)\e^{2\si l_s}|\nn |\nn P_{t-s}f|^2|^2}
{2K|\nn P_{t-s}f|^2}\Big\}(X_s)\d s.\end{split}\end{equation*}
Therefore, there exists a martingale $M_s$ for $s\in [0,t]$ such that

\beg{equation*} \beg{split} \d\eta_s^{1/2}&= \d M_s +\ff{\d \eta_s}{2\eta_s^{1/2}}-\ff{\big|
2(\U\U')(P_{t-s}f) \nn P_{t-s} f+ \ff{(\e^{2Ks}-1 ) \e^{2\si l_s}} K \nn |\nn P_{t-s}f|^2\big|^2(X_s)}
{4\eta_s^{3/2}}\\
&=\d M_s +\ff{1}{4\eta_s^{3/2}} B_s\d s,\end{split}\end{equation*}
where

\beg{equation*}\beg{split} B_s:= &2\eta_s\Big(2{ \U'}^2(P_{t-s}f) |\nn P_{t-s}f|^2 +\ff{(\e^{2Ks}-1)\e^{2\si l_s}
|\nn |\nn P_{t-s}f|^2|^2}
{2K|\nn P_{t-s}f|^2}\Big)(X_s)\\
&- \Big|
2(\U\U')(P_{t-s}f) \nn P_{t-s} f+ \ff{\e^{2Ks}-1} K \e^{2\si l_s} \nn |\nn P_{t-s}f|^2\Big|^2(X_s)\\
\ge & \ff{(\e^{2Ks}-1)\e^{2\si l_s}} K\Big\{ \ff{\U^2(P_{t-s}f)|\nn |\nn P_{t-s}f|^2|^2}
{2|\nn P_{t-s}f|^2} +4 |\nn P_{t-s}f|^4 {\U'}^2(P_{t-s}f)\\
&\qquad\qquad \qquad\qquad- 4 (\U\U')(P_{t-s}f)\<\nn P_{t-s}f, \nn |\nn P_{t-s}f|^2\>\Big\}(X_s)\\
\ge & 0.\end{split}\end{equation*} So, $\eta_s^{1/2}$ is a sub-martingale on $[0,t]$.  Therefore,
$\E\eta_0^{1/2}\le \E\eta_t^{1/2},$ which is nothing but (\ref{4.3}). \end{proof}

  \beg{thebibliography}{99}

\bibitem{B1} D. Bakry, \emph{Transformations de Riesz pour les
semigroupes sym\'etriques,} Lecture Notes in Math. No. 1123,
130--174, Springer, 1985.

\bibitem{B2} D. Bakry, \emph{On Sobolev and logarithmic
 Sobolev inequalities for
Markov semigroups,} New Trends in Stochastic Analysis, 43--75,
 World Scientific, 1997.

 \bibitem{BE}   D. Bakry and M. Emery, \emph{Hypercontractivit\'e de
semi-groupes de diffusion}, C. R. Acad. Sci. Paris. S\'er. I Math.
299(1984), 775--778.

\bibitem{BL} D. Bakry and M. Ledoux, \emph{L\'evy-Gromov's isoperimetric inequality for an infinite dimensional diffusion operator,}
Invent. Math. 123(1996), 259--281.

\bibitem{DL} H. Donnely and P. Li, \emph{Lower bounds for the
eigenvalues of Riemannian manifolds,} Michigan Math. J. 29 (1982),
149--161.

\bibitem{DS} J.-D. Deuschel and D. W. Stroock, \emph{Large
Deviations,} Academic Press, New York, 1989.

\bibitem{EL} K. D. Elworthy, \emph{Stochastic flows on
Riemannian manifolds,} Diffusion Processes and Related Problems in
Analysis, vol. II, Progress in Probability, 27, 37--72,
Birkh\"auser, 1992.

\bibitem{Hsu} E. P. Hsu, \emph{Multiplicative functional for the heat
equation on manifolds with boundary,} Michigan Math. J. 50(2002),
351--367.

\bibitem{Le} M. Ledoux, \emph{The geometry of Markov diffusion
generators,} Ann. de la Facul. des Sci. de Toulouse 9 (2000), 305--366.

\bibitem{W04} F.-Y. Wang, \emph{Equivalence of dimension-free Harnack
inequality and curvature condition,} Int. Equ. Operat. Theory 48(2004),
547--552.

\bibitem{W05} F.-Y. Wang, \emph{Gradient estimates and the first
Neumann eigenvalue on manifolds with boundary,} Stoch. Proc. Appl.
115(2005), 1475--1486.

 \bibitem{W09} F.-Y. Wang, \emph{Second fundamental form and gradient of Neumann semigroups,}  J. Funct. Anal. 256(2009), 3461--3469.

 \bibitem{W09b} F.-Y. Wang, \emph{Robin heat semigroup and    HWI inequality
on  manifolds with boundary,} preprint.

\end{thebibliography}

\end{document}